\documentclass{amsart}
\usepackage{amsmath}
\usepackage{amssymb}
\usepackage{mathrsfs}
\usepackage{graphicx}
\input xy
\xyoption{all}
\usepackage{color}
\usepackage{enumerate}
\newtheorem{theorem}{Theorem}[section]

\newtheorem{lemma}[theorem]{Lemma}

\newtheorem{proposition}[theorem]{Proposition}
\theoremstyle{definition}

\newenvironment{Proof}{{\sc Proof.}\ }{~\rule{1ex}{1ex}\vspace{0.2truecm}}

\newcommand{\Aut}{\mbox{\rm Aut}}

\newcommand{\GL}{{\rm GL}}

\begin{document}
	
	\title[Local solvability of subnormal and quasinormal subgroups of division rings]{Locally solvable subnormal and quasinormal subgroups in division rings}
	
	\author[Le Qui Danh]{Le Qui Danh}
	\address{Faculty of Mathematics and Computer Science, VNUHCM-University of Science,	227 Nguyen Van Cu Str., Dist. 5, HCM-City, Vietnam; and Department of Mathematics, Mechanics and Informatics, University of Architecture, 196 Pasteur Str., Dist. 3, HCM-City,
		Vietnam}
	\email{danh.lequi@uah.edu.vn}
	
	\author[Huynh Viet Khanh]{Huynh Viet Khanh}
	\address{Faculty of Mathematics and Computer Science, VNUHCM - University of Science, 227 Nguyen Van Cu Str., Dist. 5, Ho Chi Minh City, Vietnam.}  
	\email{huynhvietkhanh@gmail.com}

	\maketitle
\begin{abstract} 
	Let $D$ be a division ring with center $F$, and $G$ a subnormal or quasinormal subgroup of $D^*$. We show that if $G$ is locally solvable, then $G$ is contained in $F$.
\end{abstract}

\section{Introduction} 

The present paper is devoted to examining the algebraic structure of locally solvable subnormal subgroups and locally solvable quasinormal subgroups in a division ring. It shall turn out that such kinds of subgroups are always contained in the center of the division ring. 

For a moment, let us recall that a subgroup $N$ of a group $G$ is said to be \textit{subnormal} if there is a finite chain of subgroups 
$$N = N_r \leq N_{r-1}\leq\cdots\leq N_0=G,$$ for which $N_{i+1}$ is normal in $N_i$.
Whereas, if $Q$ is a subgroup of $G$ such that the relation $QH=HQ$ holds for any subgroup $H$ of $G$, then we say that $Q$ is \textit{quasinormal} (or \textit{permutable}) in $G$. 

There are close relations between the two kinds of subgroups, and we recommend to \cite[Chapter 7]{lennox-stonehewer} for additional information in the details. Here, it is noteworthy that if $G$ is finitely generated, then every quasinormal subgroup of $G$ is subnormal (\cite[Theorem B]{stonehewer}). On the other hand, we do not have the converse; there are examples demonstrating that the two notations are distinguished. To illustrate this situation, let $G$ be the dihedral group of order 8 generated by subgroups $A$ and $B$ which are of order 2. Then, it is obvious that  $AB\ne BA$ since $|AB|=4$ and $G\ne AB$, implying that $A$ and $B$ are not quasinormal subgroups of $G$. However, the nilpotency of $G$ ensures that both $A$ and $B$ are subnormal.  Furthermore, the authors in \cite{danh_19} have demonstrated that there exists a division ring which contains quasinormal subgroups that fail to be subnormal.

In the literature, there are very rich results concerning the algebraic structure multiplicative subgroups of in division ring (see \cite{hazrat}, for instance). As a direction of the study, in 1950's and 1960's, many authors paid attention on an interesting problem that to figure out how far $D^*$, and more generally its subnormal subgroups, from being abelian. In this direction, a well-known result of L. K. Hua says that if $D^*$ is solvable, then $D$ is a field.  This result of Hua has motivated many authors in attempts to examine various aspects of subnormal subgroups, in stead of $D^*$. For example, it was shown that a subnormal subgroup of $D^*$ must be contained in the center $F$ of $D$ if it is locally nilpotent, solvable, or $n$-Engel (\cite{stuth},\cite{huz},\cite{rem-kiani}), respectively). In the case of local solvability, the work of A. E. Zalesskii \cite{zalesskii} shows that every locally solvable normal subgroup of $D^*$ is contained in $F$. Now, we shall extend this result to a subnormal subgroups, rather than normal subgroups. It is unknown till now that whether every locally solvable subnormal subgroup $G$ of $D^*$ is central. A positive answer was given for only some particular cases where $D$ is assumed to be algebraic over the center (\cite{hai-thin}), or where the derived subgroup $G^{(i)}$ of $G$ is supposed to be algebraic over the center $F$ for some  $i\geq 1$ (\cite{khanh}). In section 2, we shall give the affirmative answer in the general setting to the question; that is, we show that every locally solvable subnormal subgroup is contained in $F$. And finally, in section 3, among other results, we shall figure out that the analogous fact also holds for locally solvable quasinormal subgroups. 

\section{Locally solvable subnormal subgroups}
We begin with a group-theoretic lemma which, despite its apparent simplicity, shall be frequently applied in the sequel.
\begin{lemma}\label{lemma_periodic normal}
	Every group contains a unique maximal periodic normal subgroup. Moreover, such a subgroup is characteristic in the whole group.
\end{lemma}
\begin{Proof}
	Our proof shall be obtained by mainly using Zorn's Lemma. First, we define a family of subgroups of a $G$ by taking
	$$\mathscr{A}=\{H| H \hbox{ is a periodic normal subgroup of } G\}.$$
	This family is obviously non-empty since the identity subgroup belongs to $\mathscr{A}$. Now, we consider an arbitrary chain $\{H_i\}$ of subgroups in $\mathscr{A}$. Our task, of course, is to show that $\cup\; H_i$ is again a member of $\mathscr{A}$; that is, to prove that $\cup \;H_i$ forms a periodic normal subgroup of $G$. For this purpose, pick any two elements $a, b\in \cup \;H_i$. Then, there exist indices $i$ and $j$ for which $a\in H_i$ and $b\in H_j$. Since the collection $\{H_i\}$ forms a chain, either $H_i\subseteq H_j$ or $H_j\subseteq H_i$.  It is clear that we may assume that $H_i\subseteq H_j$ and so $a.b^{-1}\in H_j\subseteq \cup \;H_i$. This implies that $\cup \;H_i$ is a subgroup of $G$. The normality as well as the periodicity of $\cup \;H_i$ may be obtained by the same way. All of this shows that $\cup \;H_i$ is a member of  $\mathscr{A}$, completing our task. Therefore, on the basic of Zorn's Lemma, the family $\mathscr{A}$ contains a maximal element $M$. 
	
	Next, we shall prove that $M$ is maximal with respect to being periodic and normal. Let $N$ be a periodic normal subgroup of $G$ for which $M\subseteq N$. Since $M$ is a maximal element of $\mathscr{A}$ and $N\in \mathscr{A}$, we must have $M=N$, which implies the maximality of $M$. 
	
	To see the uniqueness of $M$, take any periodic normal subgroup $A$ of $G$. The normality of $M$ and $A$ in $G$ permits us to form the product subgroup $AM$, which is obviously a periodic normal subgroup of $G$. But then, the maximality of $M$ reveals that $AM=M$, or $A\subseteq M$.  This argument shows every periodic normal subgroup of $G$ is contained in $M$, proving the uniqueness of $M$. 

	It only remains to show that $M$ is characteristic in $G$. For this purpose, we pick $\varphi\in \Aut(G)$, then $\varphi(M)$ is certainly a periodic normal subgroup of $G$. The uniqueness of $M$ implies that $\varphi(M)=M$. Our proof is finally finished.
\end{Proof}

For any group $G$, let us denote by $\tau(G)$ the unique maximal periodic normal subgroup of $G$ and by $B(G)$ the subgroup of $G$ such that $B(G)/\tau(G)$ is the Hirsch-Plotkin radical of $G/\tau(G)$. Phrased otherwise, $B(G)$ be the preimage of the Hirsch-Plotkin radical of the group $G/\tau(G)$ via the natural homomorphism $G\to G/\tau(G)$. 

We recall that the \textit{Hirsch-Plotkin radical} of a group is defined to be the subgroup generated by all locally nilpotent normal subgroups of the whole group. It turns out that the  Hirsch-Plotkin is the largest locally nilpotent normal subgroup. By this, the subgroup $B(G)$ is certainly normal in $G$ but, in general, not necessarily locally nilpotent. However, it is interesting to see that if $G$ is assumed to be a subnormal subgroup in a division ring, then this is the case and it is even central as we are about to see. 

\begin{proposition}\label{lemma_2.2}
	Let $D$ be a division ring with center $F$. If $G$ is a subnormal subgroup of $D^*$, then $B(G)$ is contained in $F$.
\end{proposition}

\begin{Proof}
	Being a normal subgroup of $G$, the subgroup $\tau(G)$ is a periodic subnormal subgroup of $D^*$. With reference to \cite[Theorem 8]{her}, we conclude that $\tau(G)$ is contained in $F$. 
	
	Our next step is to  assert that $B(G)$ is indeed a locally nilpotent group. For this purpose, we take an arbitrary finitely generated subgroup $H$ of $B(G)$, and our aim is to show that this is a nilpotent group. It is a simple matter to see that $H\tau(G)/\tau(G)$ is a finitely generated subgroup of $B(G)/\tau(G)$. Accordingly, the local nilpotence of $B(G)/\tau(G)$ implies that $H\tau(G)/\tau(G)$ is nilpotent. We set
	$$H_1=[H,H],\;\;\; H_2=[H_1, H] ,$$
	$$H_3=[H_{2}, H], \;\;\; \cdots\;\;\;\;\;\;\;\;\;\;\;\;\;\;\;\;$$ 
	where $[H,K]$, in particular, stands for the subgroup of $G$ generated by the set of commutators $\{[a,b]=a^{-1}b^{-1}ab|\mbox{ for all }a\in H \mbox{ and } b\in K\}$.
	Now, as $H\tau(G)/\tau(G)$ is nilpotent, we can find an integer $n$ for which $H_n\subseteq \tau(G)\subseteq F$. This fact says that any element of $H_n$ commutes elementwise with $H$ and, in consequence, we have $H_{n+1}=[H_n,H]=1$, from which it follows that $H$ is nilpotent. In other words, we obtain that $B(G)$ is locally nilpotent, as asserted.
	
	As we have pointed out before, $B(G)$ is a normal subgroup of $G$. This assures us to conclude that $B(G)$ is a locally nilpotent subnormal subgroup of $D^*$. By virtue of Huzurbazar's result, we finally obtain that $B(G)\subseteq F$. Our proof is finished.
\end{Proof}

Let $G$ be group and let $H$ be a subgroup of $G$. The normalizer of $H$ in $G$ is defined to be the subset $N_G(H)=\{g\in G|g^{-1}Hg\subseteq H\}$; and, this subset is indeed a subgroup of $G$. The following lemma, which shall provide the key to later success, gives us a way to calculate the normalizer of a locally solvable subgroup in a division ring.

\begin{lemma}[{\cite[Point 20]{wehrfritz87}}]\label{lemma_2.1}
	Let $R=F[G]$ be an algebra over the field $F$ that is a domain. If $G$ is a locally solvable, then $R$ is an Ore domain. Moreover, if we assume that $D$ is the skew field of fractions of $R$ and that $B(G)=F^*\cap G$, then $N_{D^*}(G) = GF^*$.
\end{lemma}

\begin{lemma}\label{lemma_2.3}
	Let $D$ be a division ring with center $F$. If $G$ is a locally solvable non-central subnormal subgroup of $D^*$, then $F(G)$, the division subring generated by $G$ over $F$, is coincided with $D$ and the normalizer of $G$ in $D^*$ is also locally solvable.
\end{lemma}

\begin{Proof}
	With reference to previous lemma, the local solvability of $G$ assures us to conclude that $R=F[G]$ is an Ore domain. Accordingly, its skew field of fractions is exactly $F(G)$, the division subring of $D$ generated by $G$ over $F$. Since $F(G)$ contains $G$ which is assumed to be non-central, in the light of Stuth's Theorem (\cite[Theorem 1]{stuth}), we obtain that $F(G)=D$. 
	
	Next, we argue that $B(G)=F^*\cap G$. First, it follows directly from Proposition \ref{lemma_2.2} that $B(G)\subseteq F^*\cap G$, which implies that $B(G)/\tau(G)\subseteq (F^*\cap G)/\tau(G)$. In regard to the reverse inclusion, we notice that, being the Hirsch-Plotkin radical of $G/\tau(G)$, the factor group $B(G)/\tau(G)$ is the largest locally nilpotent normal subgroup of $G/\tau(G)$. On the other hand, it is clear that $(F^*\cap G)/\tau(G)$ is an abelian normal subgroup of $G/\tau(G)$, which yields that $(F^*\cap G)/\tau(G)\subseteq B(G)/\tau(G)$. In other words, we must have $(F^*\cap G)/\tau(G)= B(G)/\tau(G)$, from which it follows that $B(G)=F^*\cap G$. Our argument is now finished. 
	
	Finally, the last assertion follows immediately from the proceeding lemma.
\end{Proof}
 Before presenting the main theorem, we need a result of Zalesskii, which can be taken from \cite{zalesskii}.
\begin{lemma}\label{lemma_2.4}
	Let $D$ be a division ring with center $F$. If $G$ is locally solvable normal subgroup of $D^*$, then $G$ is contained in $F$.
\end{lemma}

Here now is the main results of this section.

\begin{theorem}\label{theorem_2.5}
	Let $D$ be a division ring with center $F$. If $G$ is a locally solvable subnormal subgroup of $D^*$, then $G$ is contained in $F$.
\end{theorem}

\begin{Proof}
	There is nothing to be proved if $D$ is commutative. Therefore, we may suppose that $D$ is non-commutative. For purposes of contradiction, we assume that $G$ is not contained in $F$. Since $G$ is a subnormal subgroup of $D^*$, there exists a finite chain of subgroups
	$$G = G_r \leq G_{r-1}\leq\cdots\leq G_0=D^*,$$
	in which $G_i$ is normal in $G_{i-1}$ for $0\leq i \leq r$. By virtue of Lemma \ref{lemma_2.3}, we conclude that $N_{D^*}(G)$, the normalizer of $G$ in $D^*$, is a locally solvable group. The normality of $G_r$ in  $G_{r-1}$ implies that  $G_{r-1}$ is contained in $N_{D^*}(G)$ and, in consequence,  the subgroup $\leq G_{r-1}$ is locally solvable and non-central. 
	
	Repeat this procedure, now starting with $G_{r-1}$, we obtain that $G_{r-2}$ is locally solvable, too. This process must eventually terminate after finite steps, and at the final stage, we have the fact that $D^*$ is locally solvable. It follows immediately from Lemma \ref{lemma_2.4} that $D$ is commutative, which is desired contradiction. Our proof is finally completed. 
\end{Proof}

\section{Locally solvable quasinormal subgroups}
Let $G$ be a group and let $Q$ a subgroup of $G$. We say that $G$ is radical over $Q$ if for each $g$ in $G$, there is a positive integer $n$ depending upon $g$ such that $g^n$ belongs to $Q$. We prepare the way by first establishing a few results concerning on rings which are radical over subgroups.

\begin{lemma}[{\cite[Theorem 2]{gross}}]\label{lemma_gross}
	Let $Q$ be a quasinormal subgroup of $G$. If $C$ is an infinite cyclic subgroup of $G$ such that $Q\cap C=1$, then $Q$ is normal in $Q^G$ and $Q/Q_G$ is abelian.
\end{lemma}

\begin{lemma}\label{QusnSbGr} Let $G$ be a group. If $Q$ is a  quasinormal subgroup of $G$, then either $G$ is radical over $Q$ or else $Q$ is subnormal in $G$ of defect at most $2$.
\end{lemma}
\begin{Proof} To start, we assume that $G$ is not radical over $Q$. As such, we can find an element $g\in G$ for which $g^n$ fails to belong to $Q$ for every integer number $n$. Let $C$ be the cyclic subgroup of $G$ generated by the element $g$. Then, the fact that $g^n\not\in Q$ for any choice of $n$ ensures that $Q\cap C=1$. By virtue of the above lemma, we obtain that $Q$ is normal in $Q^G$, which is a normal subgroup of $G$. Phrased in another way, $Q$ is a subnormal subgroup of $G$ with the correspondent series $Q \unlhd Q^G\unlhd G$. This completes our proof. 
\end{Proof}

The next lemma, which is an interesting result of C. Faith, provides the key to establish the main result of this section.

\begin{lemma}[{\cite[Theorem B]{faith1960}}]\label{lemma_faith}
	Every division ring which is radical over a proper subring is a field.
\end{lemma}

By an analogy with C. Faith's result, a ring which is radical over a subgroup may be characterized in the following manner.

\begin{proposition}\label{DivRiRaOvSbGr} 
	Let $R$ be ring with identity $1\ne 0$ and let $G$ is a subgroup of $R^*$. If $R$ is radical over $G$, then $R$ is a division ring. 
\end{proposition}
\begin{Proof}
	To prove that $R$ is a division ring, it suffices to show that each nonzero element of $R$ is right-invertible. For this purpose, we take $x$ to be an arbitrary nonzero element of $R$. The radicality over $G$ of $x$ permits us to find an integer $n\geq1$ for which $x^n\in G$. As $G$ is a group, we can find an element $g\in G$ such that $x^ng=1$. Or, equivalently, we have $x(x^{n-1}g)=1.$ This relation shows that $x$ is right-invertible with the right inverse $x^{-1}=x^{n-1}g$. Therefore, the ring $R$ is indeed a division ring and our proposition is proved. 
\end{Proof}

\begin{lemma}\label{SbRiOfDivRiContNonAbSbGr} 
	Let $D$ be a division ring, and $G$ be a non-abelian subgroup of $D^*$. Assume that $D^*$ is radical over $G$. Then, every subring containing $G$ is coincided with $D$.
\end{lemma}
\begin{Proof} 
	For a proof by contradiction, we assume that $E$ is a proper subring of $D$ containing $G$. It is a fairly simple matter to see that $E=E[G]$. The assumption on $D^*$ assures us to deduce that $D$ is radical over $E$ and so $D$ is a field by Lemma \ref{lemma_faith}. But this contrasts to the fact that $G$ is assume to be non-abelian. 
\end{Proof}

The following theorem illustrates how the multiplicative subgroup of a division ring is affected by certain subgroups over which it is radical.

\begin{theorem}\label{DivRiRadOvLoSolSbGr} 
	Let $D$ be a division ring, and $G$ be a locally solvable subgroup of $D^*$. If $D^*$ is radical over $G$,  then $D$ is a field.
\end{theorem}
\begin{Proof} 
	Suppose, to the contrary, that $D$ is non-commutative. If $G$ is abelian, then $F(G)$ is a proper subfield over which $D$ is radical. It follows from previous lemma that $D$ is a field, which violates our supposition. We may therefore assume that $G$ is non-abelian. In the light of Lemma \ref{SbRiOfDivRiContNonAbSbGr}, we obtain that $F[G]=D$ and so $G$ possesses an abelian normal subgroup $A$ for which $G/A$ is locally finite (\cite[Point 3]{wehrfritz87.LoNilp}). This last fact ensures that $G$ is radical over $A$ and, in consequence, so is $D^*$. As a result, the division ring $D$ radical over the subfield $F(A)$, from which it follows that $D$ is a field. Again, we arrive at a desired contradiction, proving our theorem.
\end{Proof}

This may be a good place to give the main result of this section.

\begin{theorem} 
	Let $D$ be a division ring with center $F$. If $Q$ is a locally solvable quasinormal subgroup of $D^*$, then $Q$ is contained in $ F$.
\end{theorem}
\begin{Proof} 
	With reference to Lemma \ref{QusnSbGr}, we have either $Q$ is subnormal in $D^*$ or else $D^*$ is radical over $Q$. In the first event, our result follows immediately from Theorem \ref{theorem_2.5}. It remains to examine the case where $D^*$ is radical over $Q$. In this case, previous theorem says that $D$ is commutative, and our result certainly holds. Our proof is now completed.
\end{Proof}

\end{document}